\documentclass[11pt,leqno]{article}
\usepackage[active]{srcltx}

\usepackage{amsfonts,latexsym,amsmath,amssymb,amsthm}
\usepackage{fullpage}
\newtheorem{theorem}{Theorem}[section]
\newtheorem{lemma}[theorem]{Lemma}
\newtheorem{fact}[theorem]{Fact}
\newtheorem{proposition}[theorem]{Proposition}

\newtheorem{remark}[theorem]{Remark}

\theoremstyle{definition}

\theoremstyle{remark}

\newtheorem*{note*}{Note}

\numberwithin{equation}{section}

\makeatletter
\newcommand{\rank}{\mathop{\operator@font rank}}
\newcommand{\conv}{\mathop{\operator@font conv}}
\newcommand{\vol}{\mathop{\operator@font vol}}
\newcommand{\onetagright}{\tagsleft@false}
\makeatother

\newcommand{\ls}{\leqslant}
\newcommand{\gr}{\geqslant}
\renewcommand{\epsilon}{\varepsilon}

\usepackage{color}

\usepackage{hyperref}

\begin{document}
\small

\title{\bf Asymptotic shape of the convex hull of isotropic
log-concave random vectors}

\medskip

\author{Apostolos Giannopoulos, Labrini Hioni and Antonis Tsolomitis}

\date{}

\mathversion{bold}
\maketitle
\mathversion{normal}
\begin{abstract}
\footnotesize Let $x_1,\ldots ,x_N$ be independent random points distributed according to an isotrop\-ic
log-concave measure $\mu $ on ${\mathbb R}^n$, and consider the random polytope
$$K_N:={\rm conv}\{ \pm x_1,\ldots ,\pm x_N\}.$$
We provide sharp estimates for the querma\ss{}integrals
and other geometric parameters of $K_N$ in the range $cn\ls N\ls\exp (n)$; these complement previous results
from \cite{DGT1} and \cite{DGT} that were given for the range $cn\ls N\ls\exp (\sqrt{n})$. One of the basic new ingredients in
our work is a recent result of E.~Milman that determines the mean width of the centroid body $Z_q(\mu )$
of $\mu $ for all $1\ls q\ls n$.
\end{abstract}

\section{Introduction}

The purpose of this work is to add new information on the asymptotic shape of random polytopes whose
vertices have a log-concave distribution. Without loss of generality we shall assume
that this distribution is also isotropic. Recall that a convex body $K$ in ${\mathbb R}^n$
is called isotropic if it has volume $1$, it is centered, i.e.~its
center of mass is at the origin, and its inertia matrix is a multiple of the identity:
there exists a constant $L_K >0$ such that
\begin{equation}\label{isotropic-condition}\int_K\langle x,\theta\rangle^2dx =L_K^2\end{equation}
for every $\theta $ in the Euclidean unit sphere $S^{n-1}$.
More generally, a log-concave probability measure $\mu$ on ${\mathbb R}^n$ is called isotropic
if its center of mass is at the origin and its inertia matrix is the identity; in this case, the isotropic constant
of $\mu $ is defined as
\begin{equation}L_{\mu} := \sup_{x\in {\mathbb R}^n} \bigl(f_{\mu}(x)\bigr)^{1/n},\end{equation}
where $f_{\mu}$ is the density of $\mu$ with respect to the Lebesgue measure.
Note that a centered convex body $K$ of volume $1$ in ${\mathbb R}^n$ is isotropic
if and only if the log-concave probability measure $\mu_K$ with density
$x\mapsto L_K^n{\bf 1}_{K/L_K}(x)$ is isotropic.

A very well-known open question in the theory of isotropic measures is the hyperplane conjecture, which
asks if there exists an absolute constant $C>0$ such that
\begin{equation}\label{HypCon}L_n:= \sup\{ L_{\mu }:\mu \ \hbox{is an isotropic log-concave measure on}\ {\mathbb R}^n\}\ls C\end{equation}
for all $n\gr 1$. Bourgain proved in \cite{Bou} that $L_n\ls c\sqrt[4]{n}\log\! n$ (more precisely, he showed that
$L_K\ls c\sqrt[4]{n}\log\! n$ for every isotropic symmetric convex body $K$ in ${\mathbb R}^n$), while Klartag \cite{Kl}
obtained the bound $L_n\ls c\sqrt[4]{n}$. A second proof of Klartag's estimate appears in \cite{KM1}.

The study of the asymptotic shape of random polytopes whose vertices have a log-concave
distribution was initiated in \cite{DGT1} and \cite{DGT}. Given an isotropic log-concave
measure $\mu $ on ${\mathbb R}^n$, for every $N\gr n$ we consider
$N$ independent random points $x_1,\ldots ,x_N$ distributed according to $\mu $ and define the random polytope
$K_N:={\rm conv}\{ \pm x_1,\ldots ,\pm x_N\}$.
The main idea in these works was to compare $K_N$ with the $L_q$-centroid body of $\mu $
for a suitable value of $q$; roughly speaking, $K_N$ is close to the body $Z_{\log (2N/n)}(\mu )$
with high probability. Recall that the $L_q$-centroid bodies $Z_q(\mu)$, $q\gr 1$, are defined
through their support function $h_{Z_q(\mu)}$, which is given by
\begin{equation}
h_{Z_q(\mu)}(y):= \|\langle \cdot ,y\rangle \|_{L_q(\mu)} = \left(\int_{{\mathbb R}^n}|\langle x,y\rangle|^qd\mu(x)\right)^{1/q}.
\end{equation}
These bodies incorporate information about the distribution of linear functionals
with respect to $\mu$. The $L_q$-centroid bodies were introduced, under a different normalization,
by Lutwak and Zhang in \cite{Lutwak-Zhang-1997}, while in \cite{PaourisGAFA} for the first time,
and in \cite{PaourisTAMS} later on, Paouris used geometric properties
of them to acquire detailed information about the distribution of the Euclidean norm with
respect to $\mu $.

It was proved in \cite{DGT1} that, given any isotropic log-concave measure $\mu $
on ${\mathbb R}^n$ and any $cn\ls N\ls e^n$, the random polytope $K_N$ defined by $N$ independent random points
$x_1,\ldots ,x_N$ which are distributed according to $\mu $
satisfies, with high probability, the inclusion
\begin{equation}K_N\supseteq c_1Z_{\log (N/n)}(\mu )\end{equation}
(for the precise statement see Fact~\ref{fact3.2}). Then,
using the fact that the volume of the $L_q$-centroid bodies satisfies the
lower bounds $|Z_q(\mu )|^{1/n}\gr c_2\sqrt{q/n}$ if $q\ls\sqrt{n}$ and
$|Z_q(\mu )|^{1/n}\gr c_3L_{\mu }^{-1}\sqrt{q/n}$ if $\sqrt{n}\ls q\ls n$ (see Section 2),
we see that for $n\ls N\ls e^{\sqrt{n}}$ we have
\begin{equation}|K_N|^{1/n}\gr c_4\,\frac{\sqrt{\log
(2N/n)}}{\sqrt{n}},\label{eq1.6} \end{equation} while in the range
$e^{\sqrt{n}}\ls N\ls e^n$ we have
\begin{equation}|K_N|^{1/n}\gr c_5L_{\mu }^{-1}\,\frac{\sqrt{\log
(2N/n)}}{\sqrt{n}} \end{equation}with probability exponentially
close to $1$. On the other hand, one can check that for every $\alpha>1$ and $q\gr 1$,
\begin{equation}{\mathbb E}\,\big [ \sigma_n (\{\theta : h_{K_N}(\theta )\gr \alpha h_{Z_q(\mu )}(\theta
)\})\big ] \ls N\alpha^{-q},\end{equation}
where $\sigma_n$ is the rotationally invariant probability measure on the Euclidean unit
sphere $S^{n-1}$. This estimate is sufficient for some sharp upper bounds. First, for all $n\ls N\ls\exp (n)$
one has
\begin{equation}{\mathbb E}\,\big [ w(K_N)\big ]\ls c_6\,
w(Z_{\log N}(\mu )),\end{equation}where the mean width $w(C)$ of a convex body $C$ in $\mathbb R^n$ containing the origin,
is defined as twice the average of its support function
on $S^{n-1}$:
$$w(C)=\int_{S^{n-1}} h_C(\theta)\,d\sigma_n (\theta).$$
Second, one has
\begin{equation}
|K_N|^{1/n}\ls c_7\frac{\sqrt{\log (2N/n)}}{\sqrt{n}}\end{equation}
with probability greater than $1-\frac{1}{N}$, where $C>0$ is an
absolute constant.

In \cite{DGT} these results were extended to the full
family of querma\ss{}integrals $W_{n-k}(K_N)$ of $K_N$. These are
defined through Steiner's formula
\begin{equation}\label{eq1.9}
|K+tB_2^n|=\sum_{k=0}^n\binom{n}{k}W_{n-k}(K)t^{n-k},
\end{equation}
where $W_{n-k}(K)$ is the mixed volume $V(K, k; B_2^n,
n-k)$. It is more convenient to express the estimates using
a normalized variant of $W_{n-k}(K)$: for every
$1\ls k\ls n$ we set
\begin{equation}
Q_k(K)=\left (\frac{W_{n-k}(K)}{\omega_n}\right )^{1/k}=\left
(\frac{1}{\omega_k}\int_{G_{n,k}}|P_F(K)|\,d\nu_{n,k}(F)\right
)^{1/k},\label{eq1.10}
\end{equation}
where the last equality follows from Kubota's integral formula (see
Section 2 for background information on mixed volumes). Then, one has the
following results on the expectation of $Q_k(K_N)$ for all values of $k$:

\begin{theorem}[Dafnis, Giannopoulos and Tsolomitis, \cite{DGT}]\label{thmDGT}
If $n^2\ls N\ls \exp (cn)$ then for every $1\ls k\ls n$ we have
\begin{equation}
L_{\mu }^{-1}\sqrt{\log N}\lesssim\mathbb E \,\big [Q_k(K_N)\big
]\lesssim w(Z_{\log N}(K)).\label{eq:1-11-b}
\end{equation}
In the range $n^2\ls N\ls \exp (\sqrt{n})$ one has an asymptotic
formula: for every $1\ls k\ls n$,
\begin{equation}
\mathbb E\,\big [ Q_k(K_N)\big ]\simeq \sqrt{\log N}.\label{eq1.12}
\end{equation}
\end{theorem}

All these estimates remain valid for $n^{1+\delta }\ls N\ls n^2$, where $\delta\in (0,1)$ is fixed,
if we allow the constants to depend on $\delta $. Working in the
range $N\simeq n$ is possible, but requires some additional attention (see e.g. \cite{Alonso-Prochno-PAMS}
for the case of mean width).

A more careful analysis (which can be found in \cite[Theorem 1.2]{DGT})
shows that if $n^2\ls N\ls \exp (\sqrt{n})$ then, for any $s\gr 1$, a random $K_N$
satisfies, with probability greater than $1-N^{-s}$,
\begin{equation}
Q_k(K_N)\ls c_1(s)\sqrt{\log N}\label{4.1}
\end{equation}
for all $1\ls k\ls n$ and, with probability greater than $1-\exp
(-\sqrt{n})$, \begin{equation} Q_k(K_N)\gr c_8\sqrt{\log
N}\label{4.2}
\end{equation}
for all $1\ls k\ls n$, where $c_1(s)>0$ depends only on $s$, and $c_8>0$ is an absolute constant.

A natural question that arises is whether these results can be extended to the full range $cn\ls N\ls \exp (n)$
of values of $N$. If one decides to follow the approach of \cite{DGT1} and \cite{DGT} then there are two main obstacles.
The first one is that the lower bound $|Z_q(\mu )|^{1/n}\gr c\sqrt{q/n}$ is currently known only in the range $q\ls\sqrt{n}$.
In fact, proving the same for larger values of $q$ would lead to improved estimates on $L_n$ (for example, see the computation after Lemma~2.2
in \cite{KM1}). The second one was
that, until recently, a sharp estimate on the mean width of $Z_q(\mu )$ was known only for $q\ls\sqrt{n}$; G.~Paouris proved
in \cite{PaourisGAFA} that for every isotropic log-concave measure $\mu $ on $\mathbb R^n$ and any $q\ls\sqrt{n}$ one has
\begin{equation}w\bigl(Z_q(K)\bigr)\ls c_9\sqrt{q}.\end{equation}
Recently, E.~Milman \cite{EM} obtained the same upper bound (modulo
logarithmic terms) for $q$ beyond $\sqrt{n}$.

\begin{theorem}[E.~Milman, \cite{EM}]\label{thmEMilman}
For every isotropic log-concave measure $\mu $ on $\mathbb R^n$ and for all $q\in[\sqrt{n},n]$ we have
\begin{align}
w(Z_q(\mu ))\ls c_{10}\sqrt{q}\ \log^2(1+q).\label{em-zq-estimate}
\end{align}
\end{theorem}

An immediate consequence of this result is that it provides a new bound for the mean width
of an origin symmetric isotropic convex body $K$ in ${\mathbb R}^n$. In this case it is known that $Z_n(K)\supseteq cK$,
and we conclude that
\begin{align}
w(K)\ls C_1\sqrt{n}\ \log^2(1+n) L_K\label{em-estimate}
\end{align}
improving the earlier known bound $w(K)\ls C_2n^{3/4}L_K$ of Hartzoulaki, from her PhD thesis \cite{MHartzoulaki}.
We note here that not all of the logarithmic terms in (\ref{em-estimate}) can be removed,
as the example of $ B_1^n/| B_1^n|^{1/n}$ shows.

\medskip


Using E.~Milman's theorem we can show the following.

\begin{theorem}\label{thm1.3}
Let $x_1,\ldots ,x_N$ be independent random points distributed according to an isotropic log-concave
measure $\mu $ on ${\mathbb R}^n$, and consider the random polytope
$K_N:={\rm conv}\{ \pm x_1,\ldots ,\pm x_N\}$. If $\exp (\sqrt{n})\ls N\ls \exp (cn)$ then for every $1\ls k\ls n$ we
have
\begin{equation}
L_{\mu }^{-1}\sqrt{\log N}\lesssim\mathbb E \,\big [Q_k(K_N)\big ]\lesssim \sqrt{\log N} \bigl(\log\log N\bigr)^2.\label{eq1.11}
\end{equation}
\end{theorem}


Next we provide estimates for $Q_k(K_N)$ for ``most'' $K_N$:

\begin{theorem}\label{thm1.6} Let $x_1,\ldots ,x_N$ be independent random points distributed according to an isotropic log-concave
measure $\mu $ on ${\mathbb R}^n$, and consider the random polytope
$K_N:={\rm conv}\{ \pm x_1,\ldots ,\pm x_N\}$. For all $\exp (\sqrt{n}) \ls N\ls \exp(n)$ and $s\gr 1$ we have
\begin{align}
Q_k(K_N)\ls c_2(s)\sqrt{\log N}\ (\log\log N)^2\label{thm1.6:eq2},
\end{align}
for all $1\ls k<n$, with probability greater than $1-N^{-s}$.
\end{theorem}

We also provide estimates on the volume radius of a random projection $P_F(K_N)$
of $K_N$ onto $F\in G_{n,k}$ (in terms of $n,k$ and $N$) in the range $e^{\sqrt{n}}\ls N\ls e^n$; these
extend the sharp estimate ${\rm v.rad}(P_F(K_N))\simeq\sqrt{\log N}$ that was obtained in \cite{DGT} for
the case $N\ls e^{\sqrt{n}}$.

\begin{theorem}\label{th:proj-vol-radius}
If $\exp (\sqrt{n})\ls N\ls e^{cn}$ and $s\gr 1$, then a random $K_N$ satisfies with probability greater than $1-\max\{N^{-s},e^{-c_{11}\sqrt{N}}\}$
the following: for every $1\ls k\ls n$ there exists a subset $M_{n,k}$ of $G_{n,k}$ with
$\nu_{n,k}(M_{n,k})\gr 1-e^{-c_{12}k}$ such that
\begin{align}
c_{13}L_{\mu }^{-1}{\sqrt{\log N}}\ls {\rm v.rad}(P_F(K_N))&:=\left (\frac{|P_F(K_N)|}{\omega_k}\right )^{1/k}\nonumber\\ &\hspace*{5em}\ls c_3(s)\sqrt{\log N}
\bigl(\log\log N\bigr)^2\label{4.24-2}
\end{align}
for all $F\in M_{n,k}$.
\end{theorem}

In Section 4 we provide an alternative proof of an estimate of Alonso-Guti\'{e}rrez, Dafnis, Hern\'{a}ndez-Cifre and Prochno
from \cite{ADHP} on the $k$-th mean outer radius
\begin{equation}\tilde{R}_k(K_N)=\int_{G_{n,k}}R(P_F(K_N))\,d\nu_{n,k}(F)\end{equation}
of a random $K_N$, as a function of $N,n$ and $k$.

\begin{theorem}\label{thm-ADHP}Let $x_1,\ldots ,x_N$ be independent random points distributed according to an isotropic log-concave
measure $\mu $ on ${\mathbb R}^n$, and consider the random polytope
$K_N:={\rm conv}\{ \pm x_1,\ldots ,\pm x_N\}$. If $n\ls N\ls\exp (\sqrt{n})$ then, for all $1\ls k\ls n$ and $s>0$ one has
\begin{equation}\label{eq:radius}c_4(s)\max\left\{\sqrt{k},\sqrt{\log (N/n)}\right\}
\ls \tilde{R}_k(K_N)\ls c_5(s)\max\left\{\sqrt{k},\sqrt{\log N}\right\}\end{equation}
with probability greater than $1-N^{-s}$, where $c_4(s),c_5(s)$ are positive constants depending only on $s$.
\end{theorem}

We provide a formula for $\tilde{R}_k(K_N)$ which is valid for all $cn\ls N\ls\exp (n)$. This allows us to recover (and explain)
the sharp estimate of Theorem \ref{thm-ADHP} for ``small" values of $N$ and to obtain its analogue for ``large" values of $N$; see
Theorem~\ref{thm-ADHP-largeN}.

In Section 5 we obtain estimates on the regularity of the covering numbers and the dual covering numbers
of a random $K_N$. In a certain range of values of $N$, these allow us to conclude that a random $K_N$ is
in $\alpha $-regular $M$-position with $\alpha\sim 1$ (see Section 5 for definitions and terminology).

\begin{theorem}\label{regularity-covering-numbers}Let $\mu $ be an isotropic
log-concave measure on ${\mathbb R}^n$. Then, assuming that $n^2\ls N\ls \exp\bigl( (n\log n)^{2/5}\bigr)$,
we have that a random $K_N$ satisfies with probability greater than $1-N^{-1}$ the entropy estimates
$$\max\left\{\log N(K_N, tr_NB_2^n), \log N(r_{N}B_2^n, tK_N)\right\}\ls c_{14}
\frac{n(\log n)^2\log(1+t)}{t}$$
for every $t\gr 1$, where $r_N=\sqrt{\log N}$ and $c_{14}>0$ is an absolute constant.
\end{theorem}

As an application we estimate the average diameter of $k$-dimensional sections of a random $K_N$, defined by
\begin{equation}\label{eq:radius-11}\tilde{D}_k(K_N)=\int_{G_{n,k}}R(K_N\cap F)\,d\nu_{n,k}(F).\end{equation}
The discussion shows that the behavior of $\tilde{D}_k(K_N)$ is not always the same as that
of $\tilde{R}_k(K_N)$. In order to give an idea of the results, let us mention here the following simplified version.

\begin{theorem}\label{thm:average-diameter}Let $\mu $ be an isotropic
log-concave measure on ${\mathbb R}^n$ and  $a,b\in (0,1)$.
\begin{enumerate}
\item[{\rm (i)}] If $k\ls bn$ then a random $K_N$ satisfies with probability $1-N^{-1}$
$$\tilde{D}_k(K_N)\ls c_b\sqrt{\log N}\quad \textrm{if } n^2\ls N\ls \exp (\sqrt{n})$$
and
$$\tilde{D}_k(K_N)\ls c_b\sqrt{\log N}(\log\log N)^2
\quad \textrm{if } \exp (\sqrt{n})\ls N\ls \exp (n).$$
\item[{\rm (ii)}] If $k\gr an$ and $N\ls\exp((n\log n)^{2/5})$ then a random $K_N$
satisfies with probability  $1-\exp(-\sqrt{n})$
$$c_a\frac{\sqrt{\log N}}{\log^3n}\ls\tilde{D}_k(K_N).$$
\end{enumerate}
where $c_a,c_b$ are positive constants that depend only on $a$ and $b$ respectively.
\end{theorem}

We conclude this paper with a brief discussion of the interesting (open) question whether the isotropic constant of a random
$K_N$ is bounded by a constant independent from $n$ and $N$. The first class of random polytopes
$K_N$ in ${\mathbb R}^n$ for which uniform bounds were established was the class of Gaussian random polytopes.
Klartag and Kozma proved in \cite{Klartag-Kozma-2009} that if $N>n$ and if $G_1,\ldots ,G_N$ are independent standard Gaussian random vectors
in ${\mathbb R}^n$, then the isotropic constant of the random
polytope $K_N={\rm conv}\{ \pm G_1,\ldots ,\pm G_N\}$ is bounded by an absolute constant $C>0$ with probability
greater than $1-Ce^{-cn}$. The same idea works in the case where the
vertices $x_j$ of $K_N$ are distributed according to an isotropic $\psi_2$-measure $\mu $; the bound
 then depends only on the $\psi_2$-constant of $\mu $. Alonso-Guti\'{e}rrez
\cite{Alonso-2008a} and Dafnis, Giannopoulos and Gu\'{e}don \cite{Dafnis-Giannopoulos-Guedon-2010} have applied the same
more or less method to obtain a positive answer in the case where the vertices of $K_N$ are chosen from the unit sphere
or an unconditional isotropic convex body respectively. We show that, in the general isotropic log-concave case,
the method of Klartag and Kozma gives the bound $O(\sqrt{\log (2N/n)})$ if $N\ls\exp (\sqrt{n})$ (a proof along the same
lines and an extension to random perturbations of random polytopes appear in \cite{ALTJ}).

\section{Notation and background material}

We work in ${\mathbb R}^n$, which is equipped with a Euclidean
structure $\langle\cdot ,\cdot\rangle $. We denote by $\|\cdot \|_2$
the corresponding Euclidean norm, and write $B_2^n$ for the
Euclidean unit ball, and $S^{n-1}$ for the unit sphere. Volume is
denoted by $|\cdot |$. We write $\omega_n$ for the volume of $B_2^n$
and $\sigma_n $ for the rotationally invariant probability measure on
$S^{n-1}$. The Grassmann manifold $G_{n,k}$ of $k$-dimensional
subspaces of ${\mathbb R}^n$ is equipped with the Haar probability
measure $\nu_{n,k}$. Let $1\ls k\ls n$ and $F\in G_{n,k}$. We will
denote the orthogonal projection from $\mathbb R^{n}$ onto $F$ by
$P_F$. We also define $B_F=B_2^n\cap F$ and $S_F=S^{n-1}\cap F$.

The letters $c,c^{\prime }, c_1, c_2$ etc. denote absolute positive
constants whose value may change from line to line. Whenever we
write $a\simeq b$, we mean that there exist absolute constants
$c_1,c_2>0$ such that $c_1a\ls b\ls c_2a$. Similarly, if
$K,L\subseteq \mathbb R^n$ we will write $K\simeq L$ if there exist
absolute constants $c_1, c_2>0$ such that $c_{1}K\subseteq L
\subseteq c_{2}K$. We also write $\overline{A}$ for the homothetic
image of volume 1 of a convex body $A\subseteq \mathbb R^n$,
i.e.~$\overline{A}:= \frac{A}{|A|^{1/n}}$.

A convex body is a compact convex subset $C$ of ${\mathbb R}^n$
with non-empty interior. We say that $C$ is symmetric
if $-x\in C$ whenever $x\in C$. We say that $C$ is centered if it
has center of mass at the origin i.e.~$\int_C\langle
x,\theta\rangle dx=0$ for every $\theta\in S^{n-1}$. The support
function $h_C\,:\,{\mathbb R}^n\rightarrow {\mathbb R}$ of $C$ is
defined by $h_C(x )=\max\{\langle x,y\rangle :y\in C\}$. For each
$-\infty < p<\infty $, $p\neq 0$, we define the $p$-mean width of
$C$ by
\begin{equation}
w_p(C):=\left(\int_{S^{n-1}}h_{C}^p(\theta)d\sigma_n (\theta )\right)^{1/p}.\label{2.1}
\end{equation}
The mean width of $C$ is the
quantity $w(C)=w_1(C)$. The radius of $C$ is defined as
$R(C)=\max\{ \| x\|_2:x\in C\}$ and, if the origin is an interior
point of $C$, the polar body $C^{\circ }$ of $C$ is
\begin{equation}
C^{\circ }:=\{ y\in {\mathbb R}^n: \langle x,y\rangle \ls
1\;\hbox{for all}\; x\in C\}.\label{2.2}
\end{equation}
Finally, if $C$ is a symmetric convex body in $\mathbb R^n$ and $\|\cdot\|_C$ is the norm induced to ${\mathbb R}^n$
by $C$, we set $$M(C)=\int_{S^{n-1}}\|x\|_Cd\sigma_n(x)$$
and write $b(C)$ for the smallest positive constant $b$ with the property $\|x\|_C\ls b\|x\|_2$ for
all $x\in {\mathbb R}^n$. From V. Milman's proof of Dvoretzky's theorem (see \cite[Chapter 5]{AGA-book}) we know
that if $k\ls cn(M(C)/b(C))^2$ then for most $F\in G_{n,k}$ we have $C\cap F\simeq \frac{1}{M(C)}\,B_F$.

\subsection{Querma\ss{}integrals}

Let ${\cal K}_n$ denote the class of non-empty compact convex subsets of ${\mathbb R}^n$.
The relation between volume and the operations of addition and
multiplication of compact convex sets by nonnegative reals is described
by Minkowski's fundamental theorem: If $K_1,\ldots ,K_m\in
{\mathcal K}_n$, $m\in {\mathbb N}$, then the volume of
$t_1K_1+\cdots +t_mK_m$ is a homogeneous polynomial of degree $n$
in $t_i\gr 0$:
\begin{equation}
|t_1K_1+\cdots +t_mK_m|=\sum_{1\ls  i_1,\ldots ,i_n\ls  m}
V(K_{i_1},\ldots ,K_{i_n})t_{i_1}\cdots t_{i_n},\label{2.4}
\end{equation}
where the coefficients $V(K_{i_1},\ldots ,K_{i_n})$ can be chosen to
be invariant under permutations of their arguments. The
coefficient $V(K_{i_1},\ldots ,K_{i_n})$ is called the mixed
volume of the $n$-tuple $(K_{i_1},\ldots ,K_{i_n})$.

Steiner's formula is a special case of Minkowski's theorem; if $K$ is a convex body in ${\mathbb R}^n$ then the
volume of $K+tB_2^n$, $t>0$, can be expanded as a polynomial in
$t$:
\begin{equation}
|K+tB_2^n|=\sum_{k=0}^n\binom{n}{k}W_{n-k}(K)t^{n-k},\label{2.5}
\end{equation}
where $W_{n-k}(K):=V(K, k; B_2^n, n-k)$ is the $(n-k)$-th
querma\ss{}integral of $K$. It will be convenient for us to work
with a normalized variant of $W_{n-k}(K)$: for every $1\ls  k\ls
n$ we set
\begin{equation}
Q_k(K)=\left
(\frac{1}{\omega_k}\int_{G_{n,k}}|P_F(K)|\,d\nu_{n,k}(F)\right
)^{1/k}.\label{2.6}
\end{equation}
Note that $Q_1(K)=w(K)$. Kubota's integral
formula
\begin{equation}
W_{n-k}(K)=\frac{\omega_n}{\omega_k}\int_{G_{n,k}}
|P_F(K)|d\nu_{n,k}(F)\label{2.7}
\end{equation}
shows that
\begin{equation}
Q_k(K)=\left (\frac{W_{n-k}(K)}{\omega_n}\right )^{1/k}.\label{2.8}
\end{equation}
The Aleksandrov-Fenchel inequality states that if $K$, $L$,
$K_3,\ldots,K_n\in {\mathcal K}_n$, then
\begin{equation}
V(K,L,K_3,\ldots ,K_n)^2\gr V(K,K,K_3,\ldots ,K_n)
V(L,L,K_3,\ldots ,K_n).\label{2.9}
\end{equation}
This implies that the
sequence $(W_0(K),\ldots ,W_n(K))$ is log-concave: we have
\begin{equation}
W_j^{k-i}\gr W_i^{k-j}W_k^{j-i}\label{2.10}
\end{equation}
\noindent if $0\ls  i<j<k\ls  n$. Taking into account (\ref{2.8}) we
conclude that $Q_k(K)$ is a decreasing function of $k$. For the
theory of mixed volumes we refer to \cite{Schneider-book}.

\mathversion{bold}\subsection{$L_q$-centroid bodies of isotropic log-concave measures}\mathversion{normal}

We denote by ${\mathcal{P}}_n$ the class of all Borel
probability measures on $\mathbb R^n$ which are absolutely
continuous with respect to the Lebesgue measure. The density of $\mu
\in {\mathcal{P}}_n$ is denoted by $f_{\mu}$. We say that $\mu
\in {\mathcal{P}}_n$ is centered if, for all $\theta\in S^{n-1}$,
\begin{equation}
\int_{\mathbb R^n} \langle x, \theta \rangle d\mu(x) = \int_{\mathbb
R^n} \langle x, \theta \rangle f_{\mu}(x) dx = 0.
\end{equation}A measure $\mu$ on
$\mathbb R^n$ is called $\log$-concave if $\mu(\lambda
A+(1-\lambda)B) \gr \mu(A)^{\lambda}\mu(B)^{1-\lambda}$ for all compact subsets $A$
and $B$ of ${\mathbb R}^n$ and all $\lambda \in (0,1)$. A function
$f:\mathbb R^n \rightarrow [0,\infty)$ is called $\log$-concave if
its support $\{f>0\}$ is a convex set and the restriction of $\log{f}$ to it is concave.
Borell has proved in \cite{Borell-1974} that if a probability measure $\mu $ is log-concave and $\mu (H)<1$ for every
hyperplane $H$, then $\mu \in {\mathcal{P}}_n$ and its density
$f_{\mu}$ is $\log$-concave. Note that if $K$ is a convex body of volume $1$ in
$\mathbb R^n$ then the Brunn-Minkowski inequality implies that
${\bf 1}_{K} $ is the density of a $\log$-concave measure.

If $\mu $ is a $\log $-concave measure on ${\mathbb R}^n$ with density $f_{\mu}$,
we define the isotropic constant of $\mu $ by
\begin{equation}\label{definition-isotropic}
L_{\mu }:=\left (\frac{\sup_{x\in {\mathbb R}^n} f_{\mu} (x)}{\int_{{\mathbb
R}^n}f_{\mu}(x)dx}\right )^{\frac{1}{n}} [\det {\rm
Cov}(\mu)]^{\frac{1}{2n}},\end{equation} where ${\rm Cov}(\mu)$ is
the covariance matrix of $\mu$ with entries
\begin{equation}{\rm Cov}(\mu )_{ij}:=\frac{\int_{{\mathbb R}^n}x_ix_j f_{\mu}
(x)\,dx}{\int_{{\mathbb R}^n} f_{\mu} (x)\,dx}-\frac{\int_{{\mathbb
R}^n}x_i f_{\mu} (x)\,dx}{\int_{{\mathbb R}^n} f_{\mu}
(x)\,dx}\frac{\int_{{\mathbb R}^n}x_j f_{\mu}
(x)\,dx}{\int_{{\mathbb R}^n} f_{\mu} (x)\,dx}.\end{equation} Note that $L_{\mu }$
is an affine invariant of $\mu $ and does not depend on the choice of the
Euclidean structure. We say
that a $\log $-concave probability measure $\mu $ on ${\mathbb R}^n$
is isotropic if it is centered and ${\rm Cov}(\mu )$ is the identity matrix.


Recall that if $\mu$ is a log-concave probability measure on ${\mathbb R}^n$ and if $q\gr 1$ then
the $L_q$-centroid body $Z_q(\mu )$ of $\mu$ is the symmetric convex body with support
function
\begin{equation}h_{Z_q(\mu )}(y):= \left(\int_{{\mathbb R}^n} |\langle x,y\rangle|^{q}d\mu (x) \right)^{1/q}.\end{equation}
Observe that $\mu $ is isotropic if and only if it is centered and
$Z_2(\mu )=B_2^n$. From H\"{o}lder's inequality it follows that
$Z_1(\mu )\subseteq Z_p(\mu )\subseteq Z_q(\mu )$ for all $1\ls p\ls
q<\infty $. Conversely, using Borell's lemma (see \cite[Appendix III]{Milman-Schechtman-book}), one
can check that
\begin{equation}\label{reverse inclusion for Zq} Z_q(\mu )\subseteq c_1\frac{q}{p}Z_p(\mu )\end{equation}
for all $1\ls p<q$. In particular, if $\mu $ is isotropic, then $R(Z_q(\mu ))\ls c_2q$.

For any $\alpha\gr 1$ and any $\theta\in S^{n-1}$ we define the $\psi_\alpha$-norm of $x\mapsto \langle x ,\theta\rangle $ as follows:
\begin{equation} \|\langle\cdot ,\theta\rangle\|_{\psi_\alpha}:=\inf \left \{ t>0 :
\int_{{\mathbb R}^n}\exp \left(\Big (\frac{|\langle x ,\theta\rangle|}{t}\Big)^\alpha\right) \,
d\mu(x)\ls 2 \right\},
\end{equation}provided that the set on the right hand side is non-empty. We say that $\mu$ satisfies a
$\psi_\alpha$-estimate with constant $b_\alpha=b_\alpha(\theta)$ in the direction of $\theta$ if we have
\begin{equation*}
\|\langle \cdot , \theta\rangle \|_{\psi_\alpha}\ls b_\alpha \|
\langle \cdot, \theta \rangle \|_2.
\end{equation*} We say that $\mu$ is a $\psi_{\alpha }$-measure with
constant $B_\alpha>0$ if \begin{equation*}\sup_{\theta\in
S^{n-1}}\frac{\| \langle \cdot, \theta \rangle
\|_{\psi_\alpha}}{\|\langle \cdot,\theta \rangle\|_2}\ls B_\alpha.
\end{equation*}
From Borell's lemma it follows that every log-concave measure is a $\psi_1$-measure
with constant $C$, where $C$ is an absolute positive constant.

From \cite{PaourisGAFA} and \cite{PaourisTAMS}
one knows that the ``$q$-moments"
\begin{equation}
I_q(\mu):= \left(\int_{{\mathbb R}^n} \|x\|_2^qdx\right)^{1/q}, \quad q\in (-n,+\infty)\setminus\{0\},
\end{equation}
of the Euclidean norm with respect to an isotropic log-concave probability measure $\mu$ on ${\mathbb R}^n$
are equivalent to $I_2(\mu) = \sqrt{n}$ as long as $|q|\ls \sqrt{n}$. Two main consequences of this fact are: (i)
Paouris' deviation inequality
\begin{equation}\label{eq:largedeviation}\mu (\{ x\in {\mathbb R}^n:\|x\|_2\gr c_3t\sqrt{n}\})\ls \exp\left (
-t\sqrt{n}\right )\end{equation} for every $t\gr 1$, where $c_3>0$ is
an absolute constant, and (ii) Paouris' small ball probability
estimate: for any $0<\varepsilon<\varepsilon_0$, one has
\begin{equation}\label{eq:small-ball}
\mu(\{x\in \mathbb R^n : \|x\|_2<\varepsilon \sqrt{n}\})\ls
\varepsilon^{c_4\sqrt{n}},
\end{equation} where $\varepsilon_0,c_4>0$ are absolute constants.

The next theorem summarizes our knowledge on the mean width of $Z_q(\mu )$. The first statement was proved
by Paouris in \cite{PaourisGAFA}, while the second one is E.~Milman's Theorem \ref{thmEMilman}.

\begin{theorem}\label{thm2.1}
Let $\mu $ be an isotropic log-concave measure on ${\mathbb R}^n$. If $1\ls
q\ls \sqrt{n}$, then
\begin{equation}
w(Z_q(\mu ))\simeq \sqrt{q}.\label{2.20-1}
\end{equation}
Moreover, for all $q\in[\sqrt{n},n]$ we have
\begin{equation}
w(Z_q(\mu ))\ls c_5\sqrt{q}\ \log^2(1+q).\label{2.20-2}
\end{equation}
\end{theorem}

The next theorem summarizes our knowledge on the volume radius of $Z_q(\mu )$. The first statement follows from the results
of \cite{PaourisGAFA} and \cite{KM1}, while the left hand-side in the second one was obtained in
\cite{LYZ} and the right hand-side in \cite{PaourisGAFA}.

\begin{theorem}\label{thm3.3}Let $\mu $ be an isotropic log-concave measure on ${\mathbb R}^n$. If $1\ls q\ls \sqrt{n}$ then
\begin{equation}
|Z_q(\mu )|^{1/n}\simeq \sqrt{q/n},\label{3.16}
\end{equation}
while if $\sqrt{n}\ls q\ls n$ then
\begin{equation}
c_6L_{\mu }^{-1}\sqrt{q/n}\ls |Z_q(\mu )|^{1/n}\ls
c_7\sqrt{q/n}.\label{3.17}
\end{equation}
\end{theorem}

The reader may find a detailed exposition of the theory of isotropic log-concave measures in the book \cite{BGVV-book-isotropic}.

\section{Estimates for the Querma\ss{}integrals}

We start with the {\it proof of Theorem~~\upshape{\ref{thm1.3}}}. Recall that the equivalence $\mathbb E\,\big [ Q_k(K_N)\big ]\simeq \sqrt{\log N}$
in the range $n^2\ls N\ls \exp (\sqrt{n})$ was proved in \cite{DGT} (see Theorem 1.1). What is new is the right hand-side estimate in \eqref{eq1.11}.
However, in \cite{DGT} it was proved that $\mathbb E[Q_k(K_N)]\ls w(Z_{\log N}(K))$ for
the full range of $N$. So the result follows immediately by applying Theorem~\ref{thm2.1}.

\medskip

To prove Theorem~\ref{thm1.6} we will need Lemma 4.2 from \cite{DGT} which
holds true in the more general setting of isotropic log-concave random vectors.

\begin{lemma}\label{lem4.2-DGT} Let $\mu $ be an
isotropic log-concave measure on ${\mathbb R}^n$.
For every $n^2\ls N\ls\exp (cn)$ and for every $q\gr\log N$ and
$r\gr 1$, we have
\begin{equation}
\int_{S^{n-1}}\frac{h_{K_N}^q(\theta)}{h_{Z_q(\mu )}^q(\theta)}\,d\sigma_n
(\theta) \ls (c_1r)^q\label{4.4}
\end{equation} with probability greater than
$1-r^{-q}$, where $c_1>0$ is an absolute constant.
\end{lemma}

\noindent\textit{Proof of Theorem~\upshape{\ref{thm1.6}}}. Let $\exp(\sqrt{n})\ls N\ls \exp (n)$. Applying H\"older's inequality we get
\begin{align*}
w(K_N) &=\int_{S^{n-1}}h_{K_N}(\theta)\,d\sigma_n(\theta)\\
&\ls \left(\int_{S^{n-1}}\bigl(h_{Z_q(\mu )}(\theta)\bigr)^{p}\,d\sigma_n(\theta)\right)^{1/p}
\left(\int_{S^{n-1}}\left(\frac{h_{K_N}(\theta)}{h_{Z_q(\mu )}(\theta)}\right)^{q}\,d\sigma_n(\theta)\right)^{1/q}\\[1ex]
&\hspace*{2em}= w_p\bigl(Z_q(\mu )\bigr) \left(\int_{S^{n-1}}\left(\frac{h_{K_N}(\theta)}{h_{Z_q(\mu )}(\theta)}\right)^{q}\,d\sigma_n(\theta)\right)^{1/q},
\end{align*}
where $p$ is the conjugate exponent of $q$. If we now choose $q=\log N\gr\sqrt{n}$
and use Lemma~\ref{lem4.2-DGT} we arrive at
$$w(K_N) \ls c_1 r w_p\bigl(Z_q(\mu )\bigr)$$ with probability greater than
$1-r^{-q}$. Since $q=\log N$ it follows that $p<2$ and thus $w_p(Z_q(\mu ))$ is equivalent
to $w(Z_q(\mu ))$ (see \cite[Chapter 5]{AGA-book}). Using this and applying Theorem~\ref{thmEMilman} we conclude that
$$w(K_N)\ls c_2r\sqrt{\log N} \bigl(\log\log N\bigr)^2$$
with probability greater than $1-r^{-\log N}$. Choosing $r=e$ we complete the proof of
(\ref{thm1.6:eq2}). \hfill$\Box$

\medskip

We can also give estimates on the volume radius of a random projection $P_F(K_N)$
of $K_N$ onto $F\in G_{n,k}$ in terms of $n,k$ and $N$. In \cite{DGT} it was shown that if
$n^2\ls N\ls \exp (\sqrt{n})$ then, a random $K_N$
satisfies with probability greater than $1-N^{-s}$ the following:
for every $1\ls k\ls n$,
\begin{equation}
c_3\sqrt{\log N}\ls {\rm v.rad} (P_F(K_N))\ls c_4(s)\sqrt{\log
N}\label{4.24}
\end{equation}
with probability greater than $1-e^{-c_5k}$ with respect to the Haar
measure $\nu_{n,k}$ on $G_{n,k}$. We extend this result to the case $\exp (\sqrt{n})\ls N\ls\exp (n)$.

For the proof we will use Theorem 1.1 from \cite{DGT1}, which was already mentioned in the introduction.
We formulate it in the more general setting of isotropic log-concave random vectors
(the probability estimate in the statement makes use of \cite[Theorem 3.13]{ALPT}: if $\gamma >1$ and $\Gamma:\ell_2^n\to\ell_2^N$ is the random operator
$\Gamma(y)=(\langle x_1, y\rangle, \ldots \langle x_N, y\rangle)$
defined by the vertices $x_1, \ldots, x_N$ of $K_N$ then ${\mathbb
P}(\|\Gamma :\ell_2^n\to\ell_2^N\|\gr \gamma\sqrt{N})\ls\exp
(-c_0\gamma\sqrt{N})$ for all $N\gr c\gamma n$---see~\cite{DGT1} for the details).

\begin{fact}\label{fact3.2}Let $\mu $ be an isotropic log-concave measure on ${\mathbb R}^n$ and let $x_1,\ldots ,x_N$
be independent random vectors distributed according to $\mu $,
with $N\gr c_1n$ where $c_1>1$ is an absolute constant. Then, for
all $q\ls c_2\log (N/n)$ we have that
\begin{equation}K_N\supseteq c_3\,Z_q(\mu ) \label{3.2}\end{equation} with probability greater than $1-\exp
(-c_4\sqrt{N})$.
\end{fact}

\noindent \textit{Proof of Theorem~\upshape{\ref{th:proj-vol-radius}}}.
For the upper bound we use \eqref{thm1.6:eq2} and Kubota's formula to get
$$\left(\frac1{\omega_k}\int_{G_{n,k}}|P_F(K_N)|\,d\nu_{n,k} (F)\right)^{1/k}
\ls c_6(s)\sqrt{\log N} \,\bigl(\log\log N\bigr)^2 L_K.$$
Applying now Markov's inequality we get that
with probability greater than $1-t^{-k}$ with respect to the Haar measure $\nu_{n,k}$
on $G_{n,k}$ we have
$$\left(\frac{|P_F(K_N)|}{\omega_k}\right)^{1/k}\ls c_6(s)t\sqrt{\log N}
\bigl(\log\log N\bigr)^2.$$ Choosing $t=e$ proves the result.

For the lower bound integrating in polar coordinates and using H\"older's inequality
 we have
\begin{align}
\int_{G_{n,k}}\frac{|P_F^{\circ }(K_N)|}{\omega_k}\,d\nu_{n,k}(F)&=
\int_{G_{n,k}}\int_{S_F}\frac{1}{h_{P_F(K_N)}^k(\theta
)}d\sigma_F(\theta )\,d\nu_{n,k}(F) \\  &=
\int_{G_{n,k}}\int_{S_F}\frac{1}{h_{K_N}^k(\theta )}d\sigma_F(\theta
)\,d\nu_{n,k}(F)\nonumber\\ &\ls  \left
(\int_{G_{n,k}}\int_{S_F}\frac{1}{h_{K_N}^n(\theta
)}d\sigma_F(\theta )\,d\nu_{n,k}(F)\right )^{k/n}\nonumber\\  &=
\left (\int_{S^{n-1}}\frac{1}{h_{K_N}^n(\theta )}d\sigma_n (\theta
)\right )^{k/n}\nonumber\\  &= \left (\frac{|K_N^{\circ
}|}{\omega_n}\right )^{k/n}.\nonumber
\end{align}
Apply now the Blaschke-Santal\'o inequality and the fact that
$K_N\supseteq c_7Z_{\log N}(\mu )$ (with probability greater than $1-\exp (-c\sqrt{N})$ (notice that
$\log N\simeq \log N/n$ for the range of $N$ we use)) to get
\begin{equation}\left (\frac{|K_N^{\circ }|}{\omega_n}\right )^{k/n}\ls
\left (\frac{\omega_n}{|K_N|}\right )^{k/n}\ls \left
(\frac{\omega_n}{|c_7Z_{\log N}(\mu )|}\right )^{k/n}.\label{4.28}
\end{equation}
Since $\log N$ is greater than $\sqrt{n}$ we can apply the inequality
$|Z_{\log N}(K)|^{1/n}\gr cL_{\mu }^{-1}\sqrt{(\log N)/n}$ to arrive at
\begin{equation}\int_{G_{n,k}}\frac{|P_F^{\circ
}(K_N)|}{\omega_k}\,d\nu_{n,k}(F)\ls \left (\frac{c_8L_{\mu }}{\sqrt{\log
N}}\right )^k.\label{4.29}
\end{equation}
Finally, we apply Markov's inequality and the reverse Santal\'o inequality of
Bourgain and V. Milman \cite{BM} to complete
the proof.\hfill$\Box$

\section{Mean outer radii}\label{sec-radii}

For any convex body $C$ in ${\mathbb R}^n$ and any $1\ls k\ls n$, the $k$-th mean outer radius of $C$ is defined
by
\begin{equation}\label{eq:radius-1}\tilde{R}_k(C)=\int_{G_{n,k}}R(P_F(C))\,d\nu_{n,k}(F).\end{equation}
Alonso-Guti\'{e}rrez, Dafnis, Hern\'{a}ndez-Cifre and Prochno studied in \cite{ADHP} the order of growth
of $\tilde{R}_k(K_N)$ as a function of $N,n$ and $k$. Their main result is Theorem \ref{thm-ADHP}:
If $n\ls N\ls\exp (\sqrt{n})$ then, for all $1\ls k\ls n$ and $s>0$ one has
\begin{equation}\label{eq:radius-2}c_1(s)\max\left\{\sqrt{k},\sqrt{\log (N/n)}\right\}
\ls \tilde{R}_k(K_N)\ls c_2(s)\max\left\{\sqrt{k},\sqrt{\log N}\right\}\end{equation}
with probability greater than $1-N^{-s}$, where $c_1(s),c_2(s)$ are positive constants depending only on $s$.

In this section we give an alternative (and simpler) proof of this result. We also extend the estimates to the
range $\exp (\sqrt{n})\ls N\ls\exp (n)$. Our approach is based on the next general fact, which is a standard
application of concentration of measure on the Euclidean sphere (see \cite[Section 5.7]{AGA-book} for the details).
If $C$ is a symmetric convex body in $\mathbb R^n$
then, for any $1\ls k<n$ and any $s>1$ there exists a subset $\Gamma_{n,k}\subset G_{n,k}$  with
measure greater than $1-e^{-c_1s^2k}$ such that the orthogonal projection of $C$ onto any subspace $F\in \Gamma_{n,k}$ satisfies
\begin{equation}\label{eq:radius-3}
R(P_F(C)) \ls w(C)+c_2s\sqrt{k/n}R(C),
\end{equation} where $c_1>0,c_2>1$ are absolute constants. In fact, one has that the reverse inequality
$R(P_F(C))\gr c\max\{ w(C),\sqrt{k/n}R(C)\}$ holds for most $F\in G_{n,k}$.
To see this, first note that if $x\in C$ and $\|x\|_2=R(C)$ then, for most $F\in G_{n,k}$ we have $\|P_F(x)\|_2\gr c\sqrt{k/n}\|x\|_2$,
and hence $R(P_F(C))\gr c\sqrt{k/n}R(C)$; integrating with respect to $\nu_{n,k}$ we get $\tilde{R}_k(C)\gr c\sqrt{k/n}R(C)$. On the other
hand, if $\sqrt{k/n}R(C)\ls c^{\prime }w(C)$ for a small enough absolute constant $0<c^{\prime }<1$ then V. Milman's proof of Dvoretzky's theorem shows
that most $k$-dimensional projections of $C$ are isomorphic Euclidean balls of radius $w(C)$, which implies that
$\tilde{R}_k(C)\gr cw(C)$. These observations lead to the next asymptotic formula.

\begin{proposition}\label{prop:diam-rdm-proj}
Let $C$ be a symmetric convex body in ${\mathbb R}^n$. For any $1\ls k\ls n$ one has
\begin{equation}\label{eq:formula-for-radii}\tilde{R}_k(C)\simeq w(C)+\sqrt{k/n}R(C).\end{equation}
\end{proposition}

We will exploit this formula for a random $K_N$. Because of \eqref{eq:formula-for-radii} we only need to estimate $w(K_N)$ and $R(K_N)$ for a random $K_N$.
This is done in Proposition~\ref{prop:wR-smallN} and Proposition~\ref{prop:wR-largeN} below. Essential ingredients are the
deviation and small ball probability estimates \eqref{eq:largedeviation} and \eqref{eq:small-ball} of Paouris, as well
as Fact \ref{fact3.2}.

We start with the case $N\ls\exp (\sqrt{n})$.

\begin{proposition}\label{prop:wR-smallN}If $n^2\ls N\ls\exp (\sqrt{n})$ then, for any $s\gr 1$, a random $K_N$ satisfies
$$c_1\sqrt{\log N}\ls w(K_N)\ls c_2s\sqrt{\log N}$$
and
$$c_3\sqrt{n}\ls R(K_N)\ls c_4s\sqrt{n}$$
with probability greater than $1-\max\{ N^{-s},e^{-c\sqrt{n}}\}$.
\end{proposition}

\noindent {\it Proof.} In the proof of Theorem~\ref{thm1.6} we saw that, for any $n\ls N\ls\exp (n)$,
\begin{equation}\label{eq:radius-4}w(K_N) \ls c_1s w\bigl(Z_{\log N}(\mu )\bigr)\end{equation} with probability greater than
$1-N^{-s}$. Assuming that $N\ls\exp (\sqrt{n})$ we have that $\log N\ls\sqrt{n}$; then Theorem~\ref{thm2.1} and \eqref{eq:radius-4}
show that
\begin{equation}\label{eq:radius-5}w(K_N)\ls c_2s\sqrt{\log N}\end{equation}
with probability greater than $1-N^{-s}$. For the lower bound we use Fact~\ref{fact3.2}: we know that for all $N\gr c_3n$ we have
\begin{equation}K_N\supseteq c_4\,Z_{\log (N/n)}(\mu ) \label{3.2b}\end{equation}
with probability greater than $1-\exp (-c_5\sqrt{N})$. It follows that if $N\ls \exp (\sqrt{n})$ then
$$w(K_N)\gr c_4w(Z_{\log (N/n)}(\mu ))\gr c_6\sqrt{\log (N/n)}$$
with probability greater than $1-\exp (-c_7\sqrt{N})$.

For the radius of $K_N$, applying \eqref{eq:largedeviation} we see that, for any $t\gr 2$,
\begin{equation}\label{eq:radius-6}R(K_N)=\max_{1\ls j\ls N}\|x_j\|_2\ls c_8t\sqrt{n}\end{equation}
with probability greater than $1-N\exp\left ( -t\sqrt{n}\right )\geq 1-\exp (-(t-1)\sqrt{n})\gr 1-N^{-(t-1)}$.
For the lower bound, if $n^2\ls N\ls\exp (\sqrt{n})$ we use \eqref{eq:small-ball} to write
\begin{align*}{\rm Prob}(R(K_N)\ls \varepsilon_0\sqrt{n} ) &= {\rm Prob}\left (\max_{1\ls j\ls N}\|x_j\|_2\ls \varepsilon_0\sqrt{n}\right )\\
&=\left [\mu(\{x\in \mathbb R^n : \|x\|_2< \varepsilon_0\sqrt{n}\})\right ]^N \ls
e^{-c_9\sqrt{n}N},
\end{align*}
which shows that $R(K_N)\gr \varepsilon_0\sqrt{n}$ with probability greater than
$1-e^{-c_9\sqrt{n}N}$. \newline\null\hfill$\Box $

\begin{remark}\rm In fact, for the proof of the lower bound $R(K_N)\gr c\sqrt{n}$ we do not really need the
small ball probability estimate of Paouris. Lata\l a has proved in \cite{Latala} that if $\mu $ is a log-concave
probability measure on ${\mathbb R}^n$ then, for any norm $\|\cdot\|$ on ${\mathbb R}^n$
and any $0\ls t\ls 1$ one has
\begin{equation}\label{eq:latala-1}\mu (\{ x:\|x\|\ls t{\mathbb E}_{\mu }(\|x\|)\})\ls Ct,\end{equation}
where $C>0$ is an absolute constant. If we assume that $\mu $ is isotropic then we easily see that
${\mathbb E}_{\mu }(\|x\|_2)\ls\sqrt{n}$, and hence, choosing a small enough absolute constant $\varepsilon_0$ we
have by \eqref{eq:latala-1} that
\begin{equation*}\mu(\{x\in \mathbb R^n : \|x\|_2< \varepsilon_0\sqrt{n}\})\ls e^{-1}.\end{equation*}
This information is enough for our purposes.
\end{remark}

\noindent {\it Proof of Theorem~\upshape{\ref{thm-ADHP}}.} Let $N\ls \exp (\sqrt{n})$. From \eqref{eq:formula-for-radii} and
Proposition \ref{prop:wR-smallN} we get that $K_N$ satisfies with probability greater than
$1-\max\{ N^{-s},e^{-c\sqrt{n}}\}$ the following: for any $1\ls k\ls n$
\begin{align*}
\tilde{R}_k(K_N) &=\int_{G_{n,k}}R(P_F(K_N))\,d\nu_{n,k}(F)\simeq w(K_N)+\sqrt{k/n}R(K_N)\\
&\gr c_1\left (\sqrt{\log (N/n)}+\sqrt{k/n}\,\sqrt{n}\right )\simeq \max\left\{\sqrt{\log (N/n)},\sqrt{k}\right\}
\end{align*}
and similarly,
\begin{align*}
\tilde{R}_k(K_N) &=\int_{G_{n,k}}R(P_F(K_N))\,d\nu_{n,k}(F)\simeq w(K_N)+\sqrt{k/n}R(K_N)\\
&\ls c_2(s)\left (\sqrt{\log N}+\sqrt{k/n}\,\sqrt{n}\right )\ls 2c_2(s)\max\left\{\sqrt{\log N},\sqrt{k}\right\},
\end{align*}
as in \cite{ADHP}. $\hfill\Box $

\medskip

The next proposition will allow us to handle the case $\exp (\sqrt{n})\ls N\ls\exp (n)$.

\begin{proposition}\label{prop:wR-largeN}If $\exp (\sqrt{n})\ls N\ls \exp (n)$ then, for any $s\gr 1$, a random $K_N$ satisfies
$$c_1L_{\mu }^{-1}\sqrt{\log N}\ls w(K_N)\ls c_2s\sqrt{\log N}(\log\log N)^2$$
and
$$c_3\max\{\sqrt{n},R(Z_{\log N}(\mu ))\}\ls R(K_N)\ls c_3s\,\log N$$
with probability greater than $1-\max\{ N^{-s},e^{-c\sqrt{n}}\}$.
\end{proposition}

\noindent {\it Proof.} Applying again \eqref{eq:radius-4} in the range $\exp (\sqrt{n})\ls N\ls\exp (n)$ we have that
\begin{equation}\label{eq:radius-5-2}w(K_N)\ls c_2s\sqrt{\log N}(\log\log N)^2\end{equation}
from Theorem~\ref{thmEMilman}. For the lower bound we use again Fact~\ref{fact3.2}, Urysohn's inequality and \eqref{3.17}
from Theorem~\ref{thm3.3} to write
$$w(K_N)\gr c_4w(Z_{\log N}(\mu ))\gr c_4(|Z_{\log N}(\mu )|/|B_2^n|)^{1/n}\gr c_6L_{\mu }^{-1}\sqrt{\log N}$$
with probability greater than $1-\exp (-c_5\sqrt{N})$.

For the radius of $K_N$ we first use the estimate $R(K_N)\ls ct\sqrt{n}$ from \eqref{eq:radius-6}
with $t\simeq s\log N/\sqrt{n}$ to obtain the bound
$c\log N$ with probability greater than $1-N^{-s}$. For the lower bound, we show that $R(K_N)\gr c\sqrt{n}$ exactly as
in the proof of Proposition \ref{prop:wR-smallN}, and we also use the bound $R(K_N)\gr R(Z_{\log N}(\mu ))$.
$\hfill\Box $

\medskip

Using Proposition \ref{prop:wR-largeN} and Proposition \ref{prop:diam-rdm-proj} as in the proof of Theorem~\ref{thm-ADHP},
we arrive at the following estimate:

\begin{theorem}\label{thm-ADHP-largeN}Let $x_1,\ldots ,x_N$ be independent random points distributed according to an isotropic log-concave
measure $\mu $ on ${\mathbb R}^n$, and consider the random polytope
$K_N:={\rm conv}\{ \pm x_1,\ldots ,\pm x_N\}$. If $\exp (\sqrt{n})\ls N\ls\exp (n)$ then, for any $s\gr 1$ and for all $1\ls k\ls n$ one has
\begin{align*}
&c\max\left\{L_{\mu }^{-1}\sqrt{\log N},\sqrt{k},\sqrt{k/n}R(Z_{\log N}(\mu ))\right\}\\
&\hspace*{1in}\ls \tilde{R}_k(K_N)\ls Cs
\max\{\sqrt{\log N}(\log\log N)^2,\sqrt{k/n}\log N\}\end{align*}
with probability greater than $1-N^{-s}$, where $c,C>0$ are absolute constants.
\end{theorem}

In full generality one cannot expect something significantly
better: for example, if $\mu=\mu_1^n $ is the uniform measure on $B_1^n/|B_1^n|$ then $R(Z_{\log N}(\mu_1^n ))\simeq\log N$, and for large
values of $N$ (i.e. exponential in $N$) we get
$$\tilde{R}_k(K_N)\simeq \sqrt{k/n}\log N.$$
On the other hand, if $\mu $ satisfies a $\psi_2$ estimate with constant $b$ then
we know that
$L_{\mu }\ls C_1b$ (see \cite{KM1}) and
we also know that $I_n(\mu )\ls cb\sqrt{n}$ (see \cite{PaourisGAFA}),
which implies that $w(K_N)\ls R(K_N)\ls C_2b\sqrt{n}$.
Moreover, $Z_{\log N}(\mu)\subseteq b\sqrt{\log N}B_2^n$.
Thus, in this case (which e.g.~includes the case of the standard Gaussian measure) we get:

\begin{theorem}Let $x_1,\ldots ,x_N$ be independent random points distributed according to an isotropic log-concave
measure $\mu $ on ${\mathbb R}^n$ which satisfies a $\psi_2$-estimate with constant $b$, and consider the random polytope
$K_N:={\rm conv}\{ \pm x_1,\ldots ,\pm x_N\}$. If $n\ls N\ls\exp (n)$ and $s\gr 1$ then $K_N$ satisfies
with probability greater than $1-N^{-s}$
\begin{equation}\label{eq:radius-22}c_1b^{-1}\max\left\{\sqrt{k},\sqrt{\log (N/n)}\right\}
\ls \tilde{R}_k(K_N)\ls c_2(s)b\max\left\{\sqrt{k},\sqrt{\log N}\right\}\end{equation}
for all $1\ls k\ls n$, where $c_2(s)$ is a positive constant depending only on $s$.
\end{theorem}

\section{Entropy estimates and diameter of sections}\label{sec4}

For every pair of convex bodies $A$ and $B$ in $\mathbb R^n$,
the covering number $N(A,B)$ of $A$ by $B$ is defined to be the
smallest number of translates of $B$ whose union covers $A$. A fundamental theorem of
V. Milman states that there exists an absolute constant
$\beta>0$ such that every symmetric convex body $K$ in ${\mathbb R}^n$ has
a linear image $\tilde{K}$ which satisfies $|\tilde{K}|=|B_2^n|$ and
\begin{equation}\label{betaMposition}
\max\bigl\{ N(\tilde{K},B_2^n), N(B_2^n,\tilde{K}),N(\tilde{K}^{\circ
},B_2^n), N(B_2^n,\tilde{K}^{\circ })\bigr\} \ls\exp(\beta n).
\end{equation}
A convex body which satisfies the above is said to be in $M$-position with constant $\beta $.
Pisier has offered in \cite{Pisier-1989} a refined version of this result: for every $0<\alpha <2$
and every symmetric convex body $K$ in ${\mathbb R}^n$ there exists
a linear image $\tilde{K}_{\alpha }$ of $K$ such that
\begin{equation}\label{eq:four-entropy-estimates}\max\bigl\{ N(\tilde{K}_{\alpha },tB_2^n),N(B_2^n,t\tilde{K}_{\alpha }),
N(\tilde{K}_{\alpha }^{\circ },tB_2^n), N(B_2^n,t\tilde{K}_{\alpha }^{\circ })\bigr\}\ls\exp
\left (\frac{c(\alpha )n}{t^{\alpha }}\right )\end{equation} for
every $t\gr 1$, where $c(\alpha )$ depends only on $\alpha $, and
$c(\alpha )=O\big ((2-\alpha )^{-\alpha/2}\big )$ as $\alpha\to 2$.
One says that $\tilde{K}_{\alpha }$ is an $\alpha $-regular $M$-position of $K$ (we refer to
\cite[Chapter 8]{AGA-book} and \cite{Pisier-book} for a detailed exposition of these results).

In this section we will first show that if $\mu $ is an isotropic log-concave measure
on ${\mathbb R}^n$ then, for a considerably large range of values of $N$, a random
$K_N$ is in $\alpha $-regular $M$-position with $\alpha \sim 1$.
To this end, it is convenient to set $r_N=\sqrt{\log N}$: recall that if
$n^2\leq N\leq\exp(\sqrt{n})$ then ${\rm v.rad}(K_N)\simeq r_N$
for a random $K_N$ (in the case $N\gr\exp (\sqrt{n})$ one has
the weaker estimate $c_1L_{\mu }^{-1}r_N\ls {\rm v.rad}(K_N)\ls c_2r_N$).
We provide estimates for the covering numbers $N(K_N,tr_NB_2^n)$ and $N(r_NB_2^n,tK_N)$ for a random $K_N$
and for all $t\gr 1$; by the duality of entropy theorem of Artstein-Avidan, V. Milman and Szarek
\cite{Artstein-Milman-Szarek-2004}, these also determine the
covering numbers $N(r_NK_N^{\circ },tB_2^n)$ and $N(B_2^n,tr_NK_N^{\circ })$, thus completing
the proof of the four required entropy estimates in \eqref{eq:four-entropy-estimates}.

\begin{proposition}\label{regularity-covering-numbers-1}Let $\mu $ be an isotropic
log-concave measure on ${\mathbb R}^n$. Then a random $K_N$ satisfies with probability
greater than $1-N^{-1}$ the entropy estimate
$$\log N(K_N,t r_{N}B_2^n)\ls\left\{\begin{array}{ll}
    \frac{cn}{t^2} & \textrm{if }  n^2\ls N\ls \exp (\sqrt{n})\\[1ex]
    \frac{cn\log^4n}{t^2} & \textrm{if }
    \exp (\sqrt{n})\ls N\ls \exp (cn).
  \end{array}
  \right.
$$
for every $t\gr 1$, where $c>0$ is an absolute constant.
\end{proposition}

\noindent {\it Proof.} We simply recall that a random $K_N$ satisfies
$w(K_N)\ls c_1\sqrt{\log N}\simeq r_N$ for ``small" $N$, and $w(K_N)\ls c_2\sqrt{\log N}(\log\log N)^2
\simeq r_N(\log\log N)^2$ for ``large" $N$, by Proposition \ref{prop:wR-smallN} and Proposition \ref{prop:wR-largeN}
respectively. The bound for $N(K_N,tr_NB_2^n)$ is then a direct consequence of Sudakov's inequality
$$\log N(C,tB_2^n)\ls cn(w(C)/t)^2$$ which is true for every convex body $C$ in ${\mathbb R}^n$
and every $t>0$ (see e.g.\ \cite[Chapter 4]{AGA-book}). $\hfill\Box $

\medskip

We turn to estimates for the dual covering numbers $N(r_NB_2^n,tK_N)$. We will make use of the following fact
(see \cite{Giannopoulos-Stavrakakis-Tsolomitis-Vritsiou-TAMS} and \cite[Proposition 9.2.8]{BGVV-book-isotropic} or \cite{GEM-2014}
for the stronger statement below): If $\mu $ is an isotropic log-concave measure on $\mathbb R^n$, then for any $2\ls q\ls\sqrt{n}$ and for any
$1\ls t\ls \min\Big\{\sqrt{q}, c_1\frac{n\log q}{q^2}\Big\}$ we have
\begin{equation}\label{eq:dual-entropy-Zq}
\log N\bigl(\sqrt{q}B_2^n, tZ_q(\mu)\bigr)\ls c_2\frac{n(\log q)^2\log t}{t},
\end{equation}
where $c_1,c_2>0$ are absolute constants. Moreover, if $q\ls (n\log n)^{2/5}$ then \eqref{eq:dual-entropy-Zq}
holds true for all $t\gr 1$. Analogous estimates are available for larger values of $q$, but they are weaker
and do not seem to be final; so, we prefer to restrict ourselves to the next case.

\begin{proposition}\label{regularity-covering-numbers-2}Let $\mu $ be an isotropic
log-concave measure on ${\mathbb R}^n$. Then, assuming that $n^2\ls N\ls \exp\bigl((n\log n)^{2/5}\bigr)$,
we have that a random $K_N$ satisfies with probability greater than $1-\exp (-c_1\sqrt{N} )$ the entropy estimate
$$\log N(r_{N}B_2^n, tK_N)\ls c_2\frac{n(\log n)^2\log (1+t)}{t}$$
for every $t\gr 1$, where $c_1, c_2>0$ are absolute constants.
\end{proposition}

\noindent {\it Proof.} It is an immediate consequence of the fact that $K_N\supseteq c_3Z_{\log N}(\mu )$
with probability greater than $1-\exp (-c_1\sqrt{N})$. Then, we clearly have
$$\log N(r_{N}B_2^n, tK_N)\ls \log N(r_{N}B_2^n, c_3tZ_{\log N}(\mu )),$$
and the result follows from \eqref{eq:dual-entropy-Zq}. $\hfill\Box $

\medskip

\noindent \textit{Proof of Theorem~\upshape{\ref{regularity-covering-numbers}}.} By
Proposition~\ref{regularity-covering-numbers-2} $$\log N(r_N B_2^n, tK_N)
\leq c\frac{n\log^2n \log(1+t)}{t}.$$ By
Proposition~\ref{regularity-covering-numbers-1},
since $$N\leq\exp{\bigl((n\log n)^{2/5}\bigr)}\leq \exp(\theta \sqrt{n}),$$
for a suitable absolute constant $\theta>0$, we have
$$\log N(K_N, tr_N B_2^n)\leq\frac{cn}{t}$$
(here we can compensate for the extra factor $\theta$ in the exponent since
for the proof of
Proposition~\ref{regularity-covering-numbers-1} we can use the fact that
$Z_{\theta\sqrt{n}}(\mu)\subseteq \theta Z_{\sqrt{n}}(\mu)$). Combining the above bounds we get the result.\hfill$\Box$

\begin{remark}\rm Following the reasoning of \cite{DGT} one can also check that there exist absolute positive
constants $c_1$, $c_2$, $c_3$ and $c_4$ so that for every $0<t< 1$  a random $K_N$ satisfies with probability
greater than $1-N^{-1}$ the next entropy estimates:

\smallskip

\noindent \textup{(i)} If $n^2\ls N\ls \exp (\sqrt{n})$ then
\begin{align}
c_1n\log\frac{c_2}{t}\ls \log N(K_N,tr_NB_2^n)
\ls c_3n\log\frac{c_4}{t},
\end{align}
\noindent \textup{(ii)} If $\exp (\sqrt{n})\ls N\ls \exp (n)$ then
\begin{align}
c_1n\log\frac{c_2}{t}\ls \log N(K_N, t\tilde{r}_NB_2^n)
\ls c_3n\log\frac{c_4\bigl(\log\log N\bigr)^2}{t},
\end{align}
where $\tilde{r}_N:={\rm v.rad}(K_N)$ satisfies $c_5L_{\mu }^{-1}r_N\ls \tilde{r}_N\ls c_6r_N$.
\end{remark}

As an application we provide estimates for the average diameter of $k$-di\-men\-sional sections of a random $K_N$.
This parameter can be defined for any convex body $C$ in ${\mathbb R}^n$ and any $1\ls k\ls n$ as follows:
\begin{equation}\label{eq:inner-radius-11}\tilde{D}_k(C)=\int_{G_{n,k}}R(C\cap F)\,d\nu_{n,k}(F).\end{equation}
We shall use the next lemma that (in the case $\alpha =2$) can be essentially found in the article \cite{VMilman-1990a}
of V. Milman (see also \cite[Lemma 9.2.5]{BGVV-book-isotropic}):

\begin{lemma}\label{lem:milman-m*-covering}
Let $C$ be a symmetric convex body in $\mathbb R^n$ and assume that
\begin{equation}\label{eq:m*-1}
\log N(C,tB_2^n)\ls \frac{\gamma n}{t^{\alpha }}
\end{equation} for all $t\gr 1$ and some constants $\alpha >0$ and $\gamma \gr 1$.
Then, for every integer $1 \ls d < n$,
a subspace $H\in G_{n,d}$ satisfies
\begin{equation}
C\cap H^{\perp } \subseteq  c_1\alpha^{-1}\left(\frac{\gamma n}{d}\right)^{1/\alpha }
\log\left(\frac{n}{d}\right) B_{H^{\perp }}
\end{equation} with probability greater than $1-\exp (-c_2d)$,
where $c_1,c_2>0$ are absolute constants.
\end{lemma}

From Proposition \ref{regularity-covering-numbers-1} we know that a random $r_N^{-1}K_N$ satisfies the
assumption of Lemma \ref{lem:milman-m*-covering} with $\gamma\simeq 1$ if $N\ls\exp (\sqrt{n})$ and
$\gamma\simeq\log^4n$ if $N\gr\exp (\sqrt{n})$. Therefore, for any $k<n$
we have that if $N\ls\exp (\sqrt{n})$ then a $k$-dimensional section of $K_N$ has radius
\begin{equation}\label{eq:Dk-1}
R(K_N\cap F) \ls c_1\sqrt{\log N}\sqrt{\frac{n}{n-k}}
\log\left(\frac{n}{n-k}\right),
\end{equation} while if $\exp (\sqrt{n})\ls N\ls \exp (n)$ then the bound becomes
\begin{equation}\label{eq:Dk-2}
R(K_N\cap F) \ls c_1\sqrt{\log N}(\log n)^2\sqrt{\frac{n}{n-k}}
\log\left(\frac{n}{n-k}\right),
\end{equation}
both with probability greater than $1-\exp (-c_2(n-k))$,
where $c_1,c_2>0$ are absolute constants. From Proposition \ref{prop:wR-smallN} and Proposition \ref{prop:wR-largeN}
we also know that a random $K_N$ has radius
$$R(K_N)\ls c\max\{\sqrt{n},\log N\},$$
and the same bound is clearly true for all its sections $K_N\cap F$. Therefore, if $n\exp (-c_2(n-k))\ls 1$ (which is true
provided that $k< n-c_3\log n$) integration on $G_{n,k}$ shows that the bounds \eqref{eq:Dk-1} and \eqref{eq:Dk-2} hold
for $\tilde{D}_k(K_N)$ as well. Taking into account the fact that $\tilde{D}_k(K_N)\ls \tilde{R}_k(K_N)$ we conclude the following.

\begin{proposition}\label{average-diameter-1}Let $\mu $ be an isotropic
log-concave measure on ${\mathbb R}^n$. Then a random $K_N$ satisfies with probability
greater than $1-N^{-1}$ the following:

\smallskip

\noindent {\rm (i)} If $n^2\ls N\ls \exp (\sqrt{n})$ then:
\begin{enumerate}
\item If $k\ls\log N$ then $\tilde{D}_k(K_N)\ls c_1\sqrt{\log N}$.
\item If $k\gr \log N$ then $\tilde{D}_k(K_N)\ls c_1\min\left\{\sqrt{k},\sqrt{\log N}\sqrt{\frac{n}{n-k}}
\log\left(\frac{n}{n-k}\right)\right\}$,
\end{enumerate}
\noindent {\rm (ii)} If $\exp (\sqrt{n})\ls N\ls \exp (n)$ then:
\begin{enumerate}
\item If $k\ls n(\log\log N)^4/\log N$ then $\tilde{D}_k(K_N)\ls c_2\sqrt{\log N}(\log\log N)^2$.
\item If $k\gr n(\log\log N)^4/\log N$ then $$\tilde{D}_k(K_N)\ls c_2\min\left\{\sqrt{k/n}\log N,\sqrt{\log N}(\log n)^2\sqrt{\frac{n}{n-k}}
\log\left(\frac{n}{n-k}\right)\right\},$$
\end{enumerate}
where $c_1,c_2>0$ are absolute constants.
\end{proposition}

\begin{remark}\rm An alternative way to estimate the average radius of $K_N\cap F$ on $G_{n,k}$ for some values of
$k$ is given by the next theorem of Klartag and Vershynin from \cite{Klartag-Vershynin}: If $1\ls k\ls c_1n(M(C)/b(C))^2$, then
\begin{equation}\label{eq:Klartag-Vershynin}
\frac{c_2}{M(C)}\ls \left(\int_{G_{n,k}} R(C\cap F)^k\,
d\nu_{n,k}(F)\right)^{1/k} \ls \frac{c_3}{M(C)},
\end{equation} where $c_1,c_2,c_3>0$ are absolute constants.

Note that a random $K_N$ satisfies $K_N\supset Z_2(\mu )=B_2^n$ and integration in polar coordinates
combined with H\"{o}lder's inequality shows that
$$M(K_N)\gr \frac{1}{{\rm v.rad}(K_N)}\simeq \frac{1}{\sqrt{\log N}}.$$
Therefore, we may apply \eqref{eq:Klartag-Vershynin} to $K_N$: for all $1\ls k\ls cn/\log N$ we have
\begin{equation}\label{eq:radius-KV}\tilde{D}_k(K_N)\ls \left(\int_{G_{n,k}} R(C\cap F)^k\,
d\nu_{n,k}(F)\right)^{1/k} \ls \frac{c_3}{M(C)}\ls c_4\sqrt{\log N}.\end{equation}
\end{remark}

We pass now to lower bounds for $\tilde{D}_k(K_N)$. In fact, we will give a lower bound which is valid
for the radius of {\it every} section $K_N\cap F$, $F\in G_{n,k}$. We need the next lemma.

\begin{lemma}\label{lem:lower-bound-for-sections}
Let $C$ be a symmetric convex body in $\mathbb R^n$ and assume that
\begin{equation}
\log N(B_2^n,tC)\ls \frac{\gamma n}{t^{\alpha }}
\end{equation} for all $t\gr 1$ and some constants $\alpha >0$ and $\gamma \gr 1$.
Then, for every $1\ls k<n$ and any subspace $F\in G_{n,k}$ we have
\begin{equation}
R(C\cap F) \gr c\alpha \gamma^{-1/\alpha}(k/n)^{1/\alpha }.
\end{equation}where $c>0$ is an absolute constant.
\end{lemma}

\noindent {\it Proof.} Let $1\ls k<n$  and consider any $F\in G_{n,k}$.
By the duality of entropy theorem of S.~Artstein-Avidan, V.~Milman, and
S.~Szarek (see~\cite{Artstein-Milman-Szarek-2004})
the projection $P_F(C^{\circ })$ of
$C^{\circ }$ onto $F$ satisfies
\begin{equation}N(P_F(C^{\circ }),tB_F)\ls N(C^{\circ },tB_2^n)\ls
\exp (\frac{\gamma n}{k}\frac{k}{t^{\alpha }}),\end{equation}
for every $t\gr 1$. We apply Lemma \ref{lem:milman-m*-covering} for the body $P_F(C^{\circ })$ (with $\gamma^{\prime }=\gamma n/k$):
there exists $H\in G_{k,\lfloor k/2\rfloor }(F)$ such that
\begin{equation}P_F(C^{\circ })\cap H\subseteq c_1\alpha(\gamma n/k)^{1/\alpha }B_H.\end{equation}
Taking polars in $H$ we see that $P_H(C\cap F)
\supseteq c_1\alpha(k/\gamma n)^{1/\alpha }B_H$.
Using the fact that for every symmetric convex body $A$ in ${\mathbb R}^k$ and
every $H\in G_{k,s}$ we have $M(A\cap H)\ls\sqrt{k/s} M(A)$ (see \cite[Chapter 5]{AGA-book}) we get
\begin{align*}
w(C\cap F) &=M((C\cap F)^{\circ })\gr\frac{1}{\sqrt{2}}
M((C\cap F)^{\circ }\cap H)=\frac{1}{\sqrt{2}}w(P_H(C\cap F))\\
&\gr c_2\alpha(k/\gamma n)^{1/\alpha }.\end{align*}
The same lower bound holds for $R(C\cap F)$. $\hfill\Box $

\medskip

From Proposition \ref{regularity-covering-numbers-2} we know that if e.g. $n^2\ls N\ls \exp\bigl((n\log n)^{2/5}\bigr)$ then
a random $K_N$ satisfies with probability greater than $1-\exp (-c_1\sqrt{N})$ the entropy estimate
$$\log N(B_2^n, tr_N^{-1}K_N)\ls c_2\frac{n(\log n)^2\log (1+t)}{t}$$
for every $t\gr 1$, where $c_1, c_2>0$ are absolute constants.
Notice that the interesting range for $t$ is up to $n$ (otherwise $tr_N^{-1}K_N$ contains $B_2^n$)
so, we may apply Lemma \ref{lem:lower-bound-for-sections}
with $C= r_N^{-1}K_N$, $\gamma =\log^3n$ and
$\alpha =1$ to get:

\begin{proposition}\label{lower-bound-sections}Let $\mu $ be an isotropic
log-concave measure on ${\mathbb R}^n$. If $n^2\ls N\ls \exp\bigl((n\log n)^{2/5}\bigr)$
then a random $K_N$ satisfies with probability greater than $1-\exp (-c_1\sqrt{N})$
the following: for every $1\ls k<n$ and any subspace $F\in G_{n,k}$,
\begin{equation}
R(K_N\cap F) \gr c\sqrt{\log N}\frac{k}{n\log^3n},
\end{equation}
where $c>0$ is an absolute constant. The same bound holds for $\tilde{D}_k(K_N)$.
\end{proposition}

\begin{remark}\rm The question to give an upper bound for $M(K_N)$ seems open and interesting. Let us note that the analogous
question for $Z_q(\mu )$ is still open. The best known result appears in \cite{GEM-2014} (see also
\cite{Giannopoulos-Stavrakakis-Tsolomitis-Vritsiou-TAMS}): For any isotropic log-concave probability measure $\mu$
on ${\mathbb R}^n$ and any $2\ls q\ls q_0 := (n\log n)^{2/5}$ one has
\begin{equation}
M(Z_q(\mu)) \ls C\frac{\sqrt{\log q}}{\sqrt[4]{q}} .
\end{equation}
This estimate does not seem to be optimal; note that since $K_N\supseteq cZ_{\log N}(\mu )$ we also
have
\begin{equation}
M(K_N) \ls C\frac{\sqrt{\log\log N}}{\sqrt[4]{\log N}}
\end{equation}
for a random $K_N$, at least in the range $\log N\ls (n\log n)^{2/5}$.
\end{remark}

\section{Remarks on the isotropic constant}

In this last section we apply directly the method of Klartag and Kozma in order to estimate the isotropic
constant $L_{K_N}$ of a random $K_N$. The starting point is the inequality
\begin{equation}\label{eq:LK-1}|K_N|^{2/n}nL_{K_N}^2\ls \frac{1}{|K_N|}\int_{K_N}\| x\|_2^2\,dx\end{equation}
(it is well-known that this holds for any symmetric convex body in ${\mathbb R}^n$; see e.g. \cite{Milman-Pajor-1989}
or \cite[Chapter 3]{BGVV-book-isotropic}).
Assuming that $N\ls\exp (\sqrt{n})$ we know by (\ref{eq1.6}) that
\begin{equation}\label{eq:LK-2}|K_N|^{1/n}\gr c_1\frac{\sqrt{\log
(2N/n)}}{\sqrt{n}}\end{equation}with probability greater than $1-\exp (-c_2\sqrt{n})$.

We write ${\mathcal F}(K_N)$ for the family of facets of $K_N$ and we
denote by $[y_1,\ldots ,y_n]$ the convex hull of $y_1,\ldots ,y_n$.
Observe that, with probability equal to $1$, all the facets of $K_N$
are simplices and that, for all $1\ls j\ls n$, $x_j$ and $-x_j$
cannot belong to the same facet of $K_N$. Following \cite[Lemma 2.5]{Klartag-Kozma-2009} one can show
the next lemma.

\begin{lemma}\label{lem:KK}Let $F_1,\ldots ,F_M$ be the facets of $K_N$. Then,
\begin{equation}\label{eq:LK-3}\frac{1}{|K_N|}\int_{K_N}\| x\|_2^2dx\ls\frac{n}{n+2}\max_{1\ls s\ls
M}\frac{1}{|F_s|}\int_{F_s}\| u\|_2^2du.\end{equation}
\end{lemma}

Let $y_1,\ldots ,y_n\in {\mathbb R}^n$ and define $F=[y_1,\ldots
,y_n]$. Then, $F=T(\Delta^{n-1})$ where $\Delta^{n-1}=[e_1,\ldots
,e_n]$ and $T_{ij}=\langle y_j,e_i\rangle =:y_{ji}$. Assume that
$\det T\neq 0$. Then,
\begin{align*}
\frac{1}{|F|}\int_{F}\| u\|_2^2du &=
\frac{1}{|\Delta^{n-1}|}\int_{\Delta^{n-1}}\| Tu\|_2^2du\\
&= \frac{1}{|\Delta^{n-1}|}\int_{\Delta^{n-1}}\sum_{i=1}^n\left
(\sum_{j=1}^ny_{ji}u_j\right )^2\,du.
\end{align*}
Using the fact that
\begin{equation}\label{eq:LK-4}\frac{1}{|\Delta^{n-1}|}\int_{\Delta^{n-1}}(u_{j_1}u_{j_2})\,du=\frac{1+\delta_{j_1,j_2}}{n(n+1)},\end{equation} we see that \begin{equation}\label{eq:LK-5}\frac{1}{|F|}\int_{F}\| u\|_2^2du =
\frac{1}{n(n+1)}\sum_{i=1}^n\left ( \sum_{j=1}^ny_{ji}^2+\left
(\sum_{j=1}^ny_{ji}\right )^2\right ),\end{equation} from where one can conclude that
\begin{equation}\label{eq:LK-8}\frac{1}{|F|}\int_{F}\|
u\|_2^2du\ls\frac{2}{n(n+1)}\max_{\varepsilon_j=\pm
1}\left\|\varepsilon_1y_1+\cdots
+\varepsilon_ny_n\right\|_2^2.\end{equation}

Next we use a Bernstein type inequality (for a proof, see e.g. \cite[Theorem 3.5.16]{AGA-book}):

\begin{lemma}\label{lem:Bernstein}
Let $g_1,\ldots ,g_n$ be independent random variables with ${\mathbb
E}\,(g_j)=0$ on some probability space $(\Omega ,\mu )$. Assume that
$\| g_j\|_{\psi_1}\ls A$ for all $1\ls j\ls n$ and some constant $A>0$.
Then, \begin{equation}\label{eq:LK-9}{\mathbb P}\,\left\{ \left |\sum_{j=1}^na_jg_j\right |\gr
t\right\}\ls 2\exp \left (-c\min\left\{ \frac{t^2}{A^2\|
a\|_2^2},\frac{t}{A\| a\|_{\infty }}\right\}\right )\end{equation}
for every $t>0$.
\end{lemma}

We first fix $\theta\in S^{n-1}$ and a choice of signs
$\varepsilon_j=\pm 1$, and apply Lemma \ref{lem:Bernstein} to the
random variables $g_j(y_1,\ldots ,y_n)=\langle
\varepsilon_jy_j,\theta\rangle $ on $\Omega =({\mathbb R}^n,\mu )^n$. Since $\mu $ is
isotropic, we know that $\| g_j\|_{\psi_1}\ls C$. Choosing $\alpha =C_0\log (2N/n)$ we get
\begin{equation}\label{eq:LK-11}{\mathbb P}\,\left\{ \left |\langle \varepsilon_1y_1+\cdots +\varepsilon_ny_n,
\theta\rangle \right |>\alpha n\right\}\ls 2\exp (-c\alpha
n).\end{equation} Consider a $1/2$-net
${\mathcal N}$ for $S^{n-1}$ with cardinality $|{\mathcal
N}|\ls 5^n$. Then, with probability greater
than $1-\exp (-c_2\alpha n)$ we have
\begin{equation}\label{eq:LK-12}|\langle \varepsilon_1y_1+\cdots +\varepsilon_ny_n,\theta\rangle
|\ls \alpha n\end{equation} for every $\theta\in {\mathcal N}$
and every choice of signs $\varepsilon_j=\pm 1$. Using a standard
successive approximation argument, and taking into account all $2^n$
possible choices of signs $\varepsilon_j=\pm 1$, we get that, with
probability greater than $1-\exp (-c_3\alpha n)$,
\begin{equation}\label{eq:LK-13}\max_{\varepsilon_j=\pm 1}\|\varepsilon_1y_1+\cdots
+\varepsilon_ny_n\|_2 \ls C_1\alpha n.\end{equation} Now, we use the fact that
\begin{equation}\label{eq:LK-14}|{\mathcal F}(K_N)|\ls \binom{2N}{n}\ls\exp (c_3\alpha n/2)\end{equation}
provided that $C_0$ is large enough. Therefore, taking also Lemma \ref{lem:KK} and \eqref{eq:LK-8} into
account, we see that, with probability greater than $$1-|{\mathcal F}(K_N)|\exp
(-c_3\alpha n)\gr 1-\exp (-c_4\alpha n),$$ we have
\begin{equation}\label{eq:LK-15}\frac{1}{|K_N|}\int_{K_N}\| x\|_2^2dx
\ls C_2\alpha^2=C_3\log^2(2N/n),\end{equation} where $C_3>0$ is an absolute
constant. From \eqref{eq:LK-1} and \eqref{eq:LK-2} we get (with probability greater than $1-\exp (-c\sqrt{n})$)
\begin{equation}\label{eq:LK-18}L_{K_N}^2 \ls  \frac{c_4}{\log
(2N/n)}\,\frac{1}{|K_N|}\int_{K_N}\| x\|_2^2\,dx \ls
C_5\log (2N/n)\end{equation} and hence $L_{K_N}\ls C_6\sqrt{\log (2N/n)}$.

\bigskip

\noindent {\bf Acknowledgement.} The authors would like to acknowledge support from the programme ``API$\Sigma$TEIA II"
of the General Secretariat for Research and Technology of Greece.

\footnotesize
\bibliographystyle{amsplain}

\bigskip

\medskip

\noindent \textbf{MSC:} Primary 52A21; Secondary 46B07, 52A40, 60D05.

\bigskip

\bigskip

\noindent \textsc{Apostolos \ Giannopoulos}: Department of
Mathematics, University of Athens, Panepistimioupolis 157-84,
Athens, Greece.

\smallskip

\noindent \textit{E-mail:} \texttt{apgiannop@math.uoa.gr}

\bigskip

\noindent \textsc{Labrini \ Hioni}: Department of
Mathematics, University of Athens, Panepistimioupolis 157-84,
Athens, Greece.

\smallskip

\noindent \textit{E-mail:} \texttt{lamchioni@math.uoa.gr}

\bigskip

\noindent \textsc{Antonis \ Tsolomitis:} Department of Mathematics,
University of the Aegean, Karlovassi 832\,00, Samos, Greece.

\noindent {\it E-mail:} \texttt{atsol@aegean.gr}

\end{document}